\documentclass[10pt]{amsart}
\usepackage{amssymb,enumerate}
\usepackage{amsmath}
\usepackage{graphics}
\usepackage{cjmb5}
\usepackage{subfigure}
\usepackage[latin5]{inputenc}
\usepackage{amsmath}
\usepackage{graphicx}
\usepackage{subfigure}
\usepackage{caption}
\usepackage{float}
\usepackage{mathtools}
\usepackage{color}
\usepackage{cite}

\begin{document}

\title{On the numerical solution of the Klein-Gordon equation by Exponential B-spline collocation method}
\author{Ozlem Ersoy Hepson$^{a}$, Alper Korkmaz$^{b,*}$,  Idiris Dag$^{c}$ \\
$^{a}${\scriptsize Department of Mathematics \& Computer, Eskişehir Osmangazi University, 26480, Eskişehir, Turkey.}\\
$^{b,*}${\scriptsize Department of Mathematics, Çankırı Karatekin University, 18100, Çankırı, Turkey.} \\
$^{c}${\scriptsize Department of Computer Engineering, Eskişehir Osmangazi University, 26480, Eskişehir, Turkey.}}
\setcounter{page}{0}
\coordinates{}{}{}{}
\date{date}
\undate{date}
\corauthor{Alper Korkmaz; alperkorkmaz7@gmail.com}
\begin{unabstract}
In the present study, we solve initial boundary value problem construted on nonlinear Klein-Gordon equation. The collocation method on exponential cubic B-spline functions forming a set of basis for the functions defined in the same interval is set up for the numerical approach. The efficiency and validity of the proposed method are determined by computing the error between the numerical and the analytical solutions. The preservation of conserved quantities is also a good indicator of validity of the method.
\end{unabstract}

\maketitle
\pagestyle{myheadings}
\markboth{Korkmaz, Ersoy, Dag}{Exponential B-spline Collocation for the Klein-Gordon Equation}

\section{Introduction}

The Nonlinear Klein-Gordon Equation (KGE) with cubic nonlinear term is given in a generalized form
\begin{equation}
\frac{\partial ^2 u}{\partial t^2}-\frac{\partial ^2 u}{\partial x^2}+P'(u)=0 \label{KGE}
\end{equation}
where $u=u(x,t)$ and $P(u)$ is reasonable nonlinear function that is chosen as the potential energy in many cases\cite{whitham1}. The choice of $P'(u)$ can generate various equations such as Sine-Gordon, Landau-Ginzburg-Higgs, $\phi^4$ or reaction Duffing equations \cite{kim1}, arising to model many physical phenomena from quantum mechanics to wave motion. The origin of the KGE iss discussed in the Kragh's essay\cite{kragh1}. He mentions in his essay that some earlier types of the KGE was connected with the general relativity theory by siginificient physicists covering Klein, Fock and Schrödinger. The KGE is absolutely one of the most important equation that the nonlinear Schrödinger equation can be derived from this equation in the canonical form\cite{ablowitz1}. Galehouse \cite{galehouse1} derived the Klein-Gordon equation geometrically using appropriate gauge transformations. 

\noindent
The scattering theory for the more general KGE is proved by Schechter\cite{schechter1}. The same theory for the KGE is developed by following classical procedure of an equaivalent equation containing first order time in a particular finite energy normed Hilbert space\cite{weder}. Weder gave the proofs of the existence and completeness of the wave operators, the intertwining relations, and the invariance principle in that study\cite{weder}. Spectral and scattering theory of the KGE are derived by implementing related eignfunction expansions in a particular field case to give strong results like in the Schrödinger case \cite{lundberg1}. Motion of a KGE in an external field is examined by Tsukanov\cite{tsukanov1}. The linear response approximation is used to calculate the drift velocity of a kink of the KGE in study. 

\noindent
Since introduced, many types of solutions describing various physical phenomena covering solitary waves and kinks have been found for the KGE. Derivations of solutions of the KGE from the solutions of the wave equation are studied by Chambers\cite{chambers1}. Some bound state type solutions of the KGE for various potentials are set up explicitly by Fleischer and Soff\cite{fle1}. An exact formal solution for a particular nonlinear KGE is obtained by manipulating the solution of linear Klein-Gordon equation\cite{burt1}. The envelope soliton and hole solutions for the nonlinear KGE are studied by transforming it to the nonlinear Schrödinger equation\cite{sharma1}. Some solitary wave solutions of a class of field equations for systems with polynomial self-interactions reduce to plane-wave solutions of the KGE\cite{burt2}. Some bell and kink type solitary wave solutions formed by using some hyperbolic functions for the KGE are constructed with the help of the extended first kind elliptic sub-equation method\cite{huang1}. Some similar solutions have been derived by Akter and Akbar using modified simple equation method\cite{akter1}. Adomian implemented the decomposition method to the initial boundary value problems for the KGE\cite{adomian1}. The interaction of two solitons that are self-localized and weakly oscillating and damped for the particular form of the KGE are studied by Kudryavtsev\cite{kudr1}. 

\noindent
Besides some analytical solutions summarized above, various numerical methods have been applied to derive solutions of the problems related to the KGE. Strauss and Vazquez \cite{strauss1} solved some problems constructed on the KGE by the simplest central second difference method. The equation integrated numerically by using four explicit finite difference methods resulting that the energy conserving scheme is more suitable than the others to simulate the long time behaviour of the solutions\cite{jimenez1}. A collocation method based on thin plate splines radial basis functions is another proposed method to solve the KGE numerically \cite{dehghan1}. The solutions of the nonhomogenous form of the KGE are approximated by multiquadric quasi-interpolation and the integrated radial basis function network schemes\cite{sarboland1}. Classical polynomial cubic B-spline functions have also been adapted for the collocation method to derive numerical solutions of the KGE\cite{zahra1,rashidina1}.

\noindent 
In the present study, we solve some initial boundary value problems derived for the nonlinear KGE by the collocation method based on nonpolynomial B-splines, namely exponential cubic B-spline functions. 

\noindent
In order to solve the initial boundary value problem setup for the nonlinear KGE (\ref{KGE}) numerically, we first reduce the order of the derivative in time for a particular case assuming $-P'(u)=u-u^3$
\begin{equation}
\begin{array}{l}
v_{t}=u_{xx}+u-u^3 \\ 
u_{t}=v
\end{array}
\label{SYS}
\end{equation}

\noindent
To complete the usual classical mathematical statement of the problem, the initial and the boundary conditions are chosen as
\begin{equation}
\begin{aligned}
u\left( x,0\right) =f\left( x\right) ,\text{ \  \  \  \ }a\leq x\leq b \\ 
v\left( x,0\right) =g\left( x\right) ,\text{ \  \  \  \ }a\leq x\leq b
\label{IC}
\end{aligned}
\end{equation}%
and
\begin{equation}
u_{x}(a,t)=u_{x}(b,t)=0,
v_{x}(a,t)=v_{x}(b,t)=0%
\end{equation}
in the finite problem domain $[a,b]$.

\section{Exponential B-spline Collocation Method}
Let $\kappa$ be a uniformly distributed grids of the finite interval $[a,b]$, such as,
$$\kappa : x_m=a+mh, m=0,1, \ldots N$$
where $h=\frac{b-a}{N}$. Then, the exponential cubic B-splines set  
\begin{equation}
B_{m}(x)=\left\{ 
\begin{array}{lcc}
b_{2}\left ( (x_{m-2}-x)-\frac{1}{\rho} \sinh(\rho(x_{m-2}-x)) \right ) & , & [x_{m-2},x_{m-1}] \\ 
a_{1}+b_{1}(x_{m}-x)+c_{1}\exp(\rho(x_{m}-x))+d_{1}\exp(-\rho(x_{m}-x)) & , & [x_{m-1},x_{m}] \\ 
a_{1}+b_{1}(x-x_{m})+c_{1}\exp(\rho(x-x_{m}))+d_{1}\exp(-\rho(x-x_{m})) & , & [x_{m},x_{m+1}] \\ 
b_{2}\left ( (x-x_{m+2})-\frac{1}{\rho} \sinh(\rho(x-x_{m+2})) \right ) & , & [x_{m+1},x_{m+2}] \\ 
0 & , & otherwise%
\end{array}%
\right.\label{ecbs}
\end{equation}%
where $m=-1,0, \ldots N$, 
$$
\begin{array}{l}
a_{1}=\dfrac{\rho h\cosh(\rho h)}{\rho h\cosh(\rho h)-\sinh(\rho h)}, \\
b_{1}=\dfrac{\rho }{2} \dfrac{\cosh(\rho h)(\cosh(\rho h)-1)+\sinh^2(\rho h)}{(\rho h\cosh(\rho h)-\sinh(\rho h))(1-\cosh(\rho h))}, \\
b_{2}=\dfrac{\rho }{2(\rho h\cosh(\rho h)-\sinh(\rho h))}, \\
c_{1}=\dfrac{1}{4} \dfrac{\exp(-\rho h)(1-\cosh(\rho h))+\sinh(\rho h)(\exp(-\rho h)-1))}{(\rho h\cosh(\rho h)-\sinh(\rho h))(1-\cosh(\rho h))}, \\
d_{1}=\dfrac{1}{4} \dfrac{\exp(\rho h)(\cosh(\rho h)-1)+\sinh(\rho h)(\exp(\rho h)-1))}{(\rho h\cosh(\rho h)-\sinh(\rho h))(1-\cosh(\rho h))},
\end{array}
$$ 
with real parameter $\rho $ forms a basis for the functions defined in $[a,b]$ \cite{mccartin}. Each exponential cubic B-spline $B_m(x)$ has two continuous lowest derivatives defined in the interval $[x_{m-1},x_{m+2}]$. Different from polynomial B-splines \cite{korkmaz1},this basis function set has been chosen as trial functions in various methods to solve different problems arising in different engineering fields\cite{moh1,moh2,ozkdv,ozsivas,ozfish,ozreact}. 
The approximate solutions $U$ and $V$ can be written in terms of the
exponential B-splines as

\begin{equation}
\begin{aligned}
U(x,t)&=\sum_{i=-1}^{N+1}\delta _{i}B_{i}(x),\text{ } \\
V(x,t)&=%
\sum_{i=-1}^{N+1}\varphi _{i}B_{i}(x)  \label{f4}
\end{aligned}
\end{equation}%
where $\delta _{i}$  and $\varphi _{i}$ are time dependent parameters to be determined from the
complementary conditions. The nodal values of $U$, $V$ and their lowest two derivatives at the nodes can be found from (\ref{f4}) as 
\begin{eqnarray}
&&%
\begin{tabular}{l}
$U_{i}=U(x_{i},t)=\dfrac{s-\rho h}{2(\rho hc-s)}\delta _{i-1}+\delta _{i}+\dfrac{s-\rho h%
}{2(\rho hc-s)}\delta _{i+1},$ \\ 
$U_{i}^{\prime }=U^{\prime }(x_{i},t)=\dfrac{\rho (1-c)}{2(\rho hc-s)}\delta _{i-1}+%
\dfrac{\rho (c-1)}{2(\rho hc-s)}\delta _{i+1}$ \\ 
$U_{i}^{\prime \prime }=U^{\prime \prime }(x_{i},t)=\dfrac{\rho ^{2}s}{2(\rho hc-s)}%
\delta _{i-1}-\dfrac{\rho ^{2}s}{\rho hc-s}\delta _{i}+\dfrac{\rho ^{2}s}{2(\rho hc-s)}%
\delta _{i+1}.$%
\end{tabular}
\label{f5} \\
&&%
\begin{tabular}{l}
$V_{i}=V(x_{i},t)=\dfrac{s-\rho h}{2(\rho hc-s)}\phi _{i-1}+\phi _{i}+\dfrac{s-\rho h}{%
2(\rho hc-s)}\phi _{i+1},$ \\ 
$V_{i}^{\prime }=V^{\prime }(x_{i},t)=\dfrac{\rho (1-c)}{2(\rho hc-s)}\phi _{i-1}+%
\dfrac{\rho (c-1)}{2(\rho hc-s)}\phi _{i+1}$ \\ 
$V_{i}^{\prime \prime }=V^{\prime \prime }(x_{i},t)=\dfrac{\rho ^{2}s}{2(\rho hc-s)}%
\phi _{i-1}-\dfrac{\rho ^{2}s}{\rho hc-s}\phi _{i}+\dfrac{\rho ^{2}s}{2(\rho hc-s)}\phi
_{i+1}.$
\end{tabular}
\label{f6}
\end{eqnarray}
where $c=\cosh{\rho h}$ and $s=\sinh{\rho h}$.

\noindent
The time discretization of the system (\ref{SYS}) by usual forward finite
difference and the Crank-Nicolson method leads to
\begin{equation}
\begin{array}{r}
\dfrac{V^{n+1}-V^{n}}{\Delta t}-\dfrac{U_{xx}^{n+1}+U_{xx}^{n}}{2}-\dfrac{%
U^{n+1}+U^{n}}{2}+\dfrac{(U^{3})^{n+1}+(U^{3})^{n}}{2}=0 \\ 
\dfrac{U^{n+1}-U^{n}}{\Delta t}-\dfrac{V^{n+1}+V^{n}}{2}=0%
\end{array}
\label{e5}
\end{equation}%
where $U^{n+1}=U(x,(n+1)\Delta t)$ represents the solution at the $(n+1)$.th
time level. Linearizing the term $(U^{3})^{n+1}$ in (\ref{e5}) as
\begin{equation}
(U^{3})^{n+1}=3U^{n+1}(U^{2})^{n}-2(U^{3})^{n}  \label{e6}
\end{equation}%
yields the time-integrated system:%
\begin{equation}
\begin{array}{r}
\dfrac{V^{n+1}-V^{n}}{\Delta t}-\dfrac{U_{xx}^{n+1}+U_{xx}^{n}}{2}-\dfrac{%
U^{n+1}+U^{n}}{2}+\dfrac{3U^{n+1}(U^{2})^{n}-(U^{3})^{n}}{2}=0 \\ 
\dfrac{U^{n+1}-U^{n}}{\Delta t}-\dfrac{V^{n+1}+V^{n}}{2}=0%
\end{array}
\label{e7}
\end{equation}%

\noindent
Approximating to $U$ and $V$ by using the nodal values (\ref{f5}) and (\ref{f6}) gives
\begin{eqnarray}
&&\nu _{m1}\delta _{m-1}^{n+1}+\nu _{m2}\phi _{m-1}^{n+1}+\nu _{m3}\delta
_{m}^{n+1}+\nu _{m4}\phi _{m}^{n+1}+\nu _{m1}\delta _{m+1}^{n+1}+\nu
_{m2}\phi _{m+1}^{n+1}  \label{e9} \\
&=&\nu _{m5}\delta _{m-1}^{n}+\nu _{m2}\phi _{m-1}^{n}+\nu _{m6}\delta
_{m}^{n}+\nu _{m4}\phi _{m}^{n}+\nu _{m5}\delta _{m+1}^{n}+\nu _{m2}\phi
_{m+1}^{n}  \notag
\end{eqnarray}%
\begin{eqnarray}
&&\nu _{m2}\delta _{m-1}^{n+1}+\nu _{m7}\phi _{m-1}^{n+1}+\nu _{m4}\delta
_{m}^{n+1}+\nu _{m8}\phi _{m}^{n+1}+\nu _{m2}\delta _{m+1}^{n+1}+\nu
_{m7}\phi _{m+1}^{n+1}  \label{e10} \\
&=&\nu _{m2}\delta _{m-1}^{n}-\nu _{m7}\phi _{m-1}^{n}+\nu _{m4}\delta
_{m}^{n}-\nu _{m8}\phi _{m}^{n}+\nu _{m2}\delta _{m+1}^{n}-\nu _{m7}\phi
_{m+1}^{n}  \notag
\end{eqnarray}

\noindent
The coefficients of iteration system (\ref{e9})-(\ref{e10}) are
\begin{equation*}
\begin{array}{l}
\nu _{m1}=\left( 3K^{2}-1\right) \alpha _{1}-\gamma _{1} \\ 
\nu _{m2}=\dfrac{2}{\Delta t}\alpha _{1} \\ 
\nu _{m3}=\left( 3K^{2}-1\right) \alpha _{2}-\gamma _{2} \\ 
\nu _{m4}=\dfrac{2}{\Delta t}\alpha _{2} \\ 
\nu _{m5}=\left( 1+K^{2}\right) \alpha _{1}+\gamma _{1} \\ 
\nu _{m6}=\left( 1+K^{2}\right) \alpha _{2}+\gamma _{2} \\ 
\nu _{m7}=-\alpha _{1} \\ 
\nu _{m8}=-\alpha _{2}%
\end{array}%
\end{equation*}

where%
\begin{equation*}
K=\alpha _{1}\delta _{i-1}+\alpha _{2}\delta _{i}+\alpha _{1}\delta _{i+1}
\end{equation*}

and%
\begin{eqnarray*}
\alpha _{1} &=&\dfrac{s-ph}{2(phc-s)},\text{ }\alpha _{2}=1, \\
\gamma _{1} &=&\dfrac{p^{2}s}{2(phc-s)},\text{ }\gamma _{2}=-\dfrac{p^{2}s}{%
phc-s}
\end{eqnarray*}%

\noindent
The system system (\ref{e9})-(\ref{e10}) can be written in the matrix form
\begin{equation}
\mathbf{Ax}^{n+1}=\mathbf{Bx}^{n}  \label{e11}
\end{equation}%
where%
\begin{equation*}
\mathbf{A=}%
\begin{bmatrix}
\nu _{m1} & \nu _{m2} & \nu _{m3} & \nu _{m4} & \nu _{m1} & \nu _{m2} &  & 
&  &  \\ 
\nu _{m2} & \nu _{m7} & \nu _{m4} & \nu _{m8} & \nu _{m2} & \nu _{m7} &  & 
&  &  \\ 
&  & \nu _{m1} & \nu _{m2} & \nu _{m3} & \nu _{m4} & \nu _{m1} & \nu _{m2} & 
&  \\ 
&  & \nu _{m2} & \nu _{m7} & \nu _{m4} & \nu _{m8} & \nu _{m2} & \nu _{m7} & 
&  \\ 
&  &  & \ddots & \ddots & \ddots & \ddots & \ddots & \ddots &  \\ 
&  &  &  & \nu _{m1} & \nu _{m2} & \nu _{m3} & \nu _{m4} & \nu _{m1} & \nu
_{m2} \\ 
&  &  &  & \nu _{m2} & \nu _{m7} & \nu _{m4} & \nu _{m8} & \nu _{m2} & \nu
_{m7}%
\end{bmatrix}%
\end{equation*}%
and%
\begin{equation*}
\mathbf{B=}%
\begin{bmatrix}
\nu _{m5} & \nu _{m2} & \nu _{m6} & \nu _{m4} & \nu _{m5} & \nu _{m2} &  & 
&  &  \\ 
\nu _{m2} & -\nu _{m7} & \nu _{m4} & -\nu _{m8} & \nu _{m2} & -\nu _{m7} & 
&  &  &  \\ 
&  & \nu _{m5} & \nu _{m2} & \nu _{m6} & \nu _{m4} & \nu _{m5} & \nu _{m2} & 
&  \\ 
&  & \nu _{m2} & -\nu _{m7} & \nu _{m4} & -\nu _{m8} & \nu _{m2} & -\nu _{m7}
&  &  \\ 
&  &  & \ddots & \ddots & \ddots & \ddots & \ddots & \ddots &  \\ 
&  &  &  & \nu _{m5} & \nu _{m2} & \nu _{m6} & \nu _{m4} & \nu _{m5} & \nu
_{m2} \\ 
&  &  &  & \nu _{m2} & -\nu _{m7} & \nu _{m4} & -\nu _{m8} & \nu _{m2} & 
-\nu _{m7}%
\end{bmatrix}%
\end{equation*}

\noindent
This system consists of $2N+2$ linear equations with $2N+6$ unknown parameters $%
\mathbf{x}^{n+1}=(\delta _{-1}^{n+1},\phi _{-1}^{n+1},\delta _{0}^{n+1},\phi
_{0}^{n+1}\ldots ,\delta _{N+1}^{n+1},\phi _{N+1}^{n+1})$. A unique solution
of the system can be obtained by imposing the boundary conditions $%
U_{x}(a,t)=0,U_{x}(b,t)=0,V_{x}(a,t)=0,V_{x}(b,t)=0$ to have the following
the equations:%
\begin{equation*}
\delta _{-1}=\delta _{1},\text{ }\phi _{-1}=\phi _{1},\text{ }\delta
_{N-1}=\delta _{N+1},\text{ }\phi _{N-1}=\phi _{N+1}
\end{equation*}

\noindent
Elimination of the parameters $\delta _{-1},\phi _{-1},\delta _{N+1},\phi
_{N+1},$ using the Eqs,(\ref{e9}) and (\ref{e10}) from the the system gives
a solvable system of $2N+2$ linear equation including $2N+2$ unknown
parameters. Placing solution parameters when is computed from the system via
a variant of the Thomas algorithm gives the approximate solution over the
subregion $[x_{i},x_{i+1}].$ We need the initial parameter vectors $%
d_{1}=(\delta _{-1},\delta _{0},..\delta _{N},\delta _{N+1})$, $d_{2}=(\phi
_{-1},\phi _{0},..\phi _{N},\phi _{N+1})$ to start the iteration process for
the system of equation: To do that, the following requirements help to
determine initial parameters:

$%
\begin{array}{l}
U_x(a,0)=0=\delta _{-1}^{0}-\delta _{1}^{0}, \\ 
U_{x}(x_{i},0)=\delta _{i-1}^{0}-\delta _{i+1}^{0},i=1,...,N-1\\ 
U_{x}(b,0)=0=\delta _{N-1}^{0}-\delta _{N+1}^{0}, \\ 
V_{x}(a,0)=0=\phi _{-1}^{0}-\phi _{1}^{0} \\ 
V_{x}(x_{i},0)=\phi _{i-1}^{0}-\phi _{i+1}^{0},i=1,...,N-1 \\ 
V_{x}(b,0)=0=\phi _{N-1}^{0}-\phi _{N+1}^{0}%
\end{array}%
$
\section{Numerical Example}
We illustrate the performance of the proposed method in terms of the
accuracy. The absolute error between the obtained numerical solution and the exact solution is computed for different values of time step length $\Delta t$ and step size $h$.
The discrete maximum norm 
\begin{equation*}
L_{\infty }=\left \vert w-W\right \vert _{\infty }=\max \limits_{j}\left
\vert w_{j}-W_{j}\right \vert
\end{equation*}%
is computed where $w$ and $W$ represent the analytical and numerical solutions, respectively. 

\noindent
The conservation laws defined as \cite{zaki1}
\begin{equation*}
\begin{array}{l}
E=\frac{1}{2}\int \limits_{a}^{b}\left[ V^{2}+U_{x}^{2}-U^{2}+\frac{1}{2}%
U^{4}\right] dx \\ 
P=\frac{1}{2}\int \limits_{a}^{b}VU_{x}dx%
\end{array}%
\end{equation*}
are also computed for each time step to observe if they are preserved.

\begin{figure}[htp]
    \subfigure[$\rho =1$]{
   \includegraphics[scale =0.4] {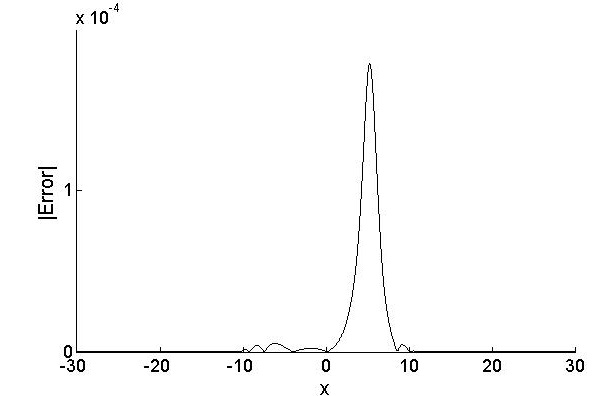}
   \label{fig:1}
 }
 \subfigure[$\rho =0.0000589$]{
   \includegraphics[scale =0.4] {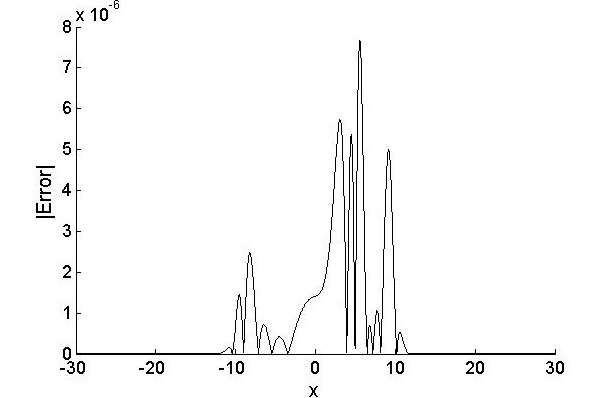}
   \label{fig:2}
 }
 \caption{Comparison of error norms $L_{\infty}$ for various values of $\rho$}
\end{figure}

\noindent
The KGE equation has an analytical solution representing a travelling wave of the form
\begin{equation}
u(x,t)=\tanh (\frac{(x-ct)}{\sqrt{2(1-c^{2})}})  \label{n1}
\end{equation}%
where $c$ is the velocity of wave\cite{zaki1}. The initial condition is generated by substituting $t=0$ into the analytical solution. In this test problem, the parameters are chosen as $\Delta t=0.005,h=0.02$, $%
c=0.5$ and various values of the exponential B-spline parameter $\rho$ in the interval $[a,b]=[-30,30]$ for  $t \in [0,30]$. The maximum error comparison at the simulation terminating time $t=30$ shows that the results are better in two decimal digits when $\rho =0.0000589$, Fig \ref{fig:1}-\ref{fig:2}.
The maximum errors are tabulated for various values of space and time step lengths and the exponential B-spline parameter $\rho$ at time $t=10$ in Table 1. A comparison with Boz's work, based on polynomial B-spline collocation, shows that the exponential B-spline generates the results one decimal digit better even when the same parameters are used.

\begin{tabular}{lllll}
\multicolumn{5}{l}{\textbf{Table 1: }$L_{\infty }\times 10^{3}$ error norms for $%
c=0.5,$ $t=10,$ $-30\leq x\leq 30$} \\ 
\multicolumn{5}{l}{} \\ \hline
Time &  & Present $(\rho=1)$ & Present $($various $\rho)$ & $\text{Ref: \cite{boz}}
$ \\ \hline
$h$ & $\Delta t$ &  &  &  \\ 
$0.2$ & $0.05$ & $17.5000$ & $3.6155(\rho=0.0000012)$ & $11.0155$ \\ 
$0.1$ & $0.02$ & $4.5179$ & $0.2525(\rho=0.0000053)$ & $2.8275$ \\ 
$0.05$ & $0.01$ & $1.1333$ & $0.0583(\rho=0.0000151)$ & $0.6611$ \\ 
$0.02$ & $0.005$ & $0.1790$ & $0.0076(\rho=0.0000589)$ & $0.0466$ \\ \hline
\end{tabular}%

\noindent
The analytical values of the conservation laws initially are calculated as $E=-13.91133798$ and $P=-0.5443310539$. The relative changes of the conservation laws defined as
\begin{equation}
\begin{aligned}
C_P(t)&=\left| \dfrac{P_{t}-P_0}{P_0}\right| \\
C_E(t)&=\left| \dfrac{E_{t}-E_0}{E_0}\right| \\
\end{aligned}
\end{equation}
where $C_P(t)$ and $C_E(t)$ are the relative changes of conservation laws at the time $t$ with respect to the reference initial values $E_0$ and $P_0$. We record the computed values of relative changes in Table 2. Even when the results obtained various values when $\rho=1$ can be recovered by changing it, Table 1, a significant recovery is not observed in the relative changes of conservation laws, Table 2.

\begin{tabular}{llllll}
\multicolumn{6}{l}{\textbf{Table 2: } Relative changes of conservative laws for various values of $\rho$} \\ 
\multicolumn{6}{l}{} \\ \hline
Time &  & $(\rho=1)$ & $(\rho=1)$ & $($various $\rho)$ & $($various $\rho)$ \\ \hline
$h$ & $\Delta t$ & $C_P(10)$ & $C_E(10)$ & $C_P(10)$ & $C_E(10)$ \\ 
$0.2$ & $0.05$ & $3.9512\times 10^{-5}$ & $1.3011\times 10^{-6}$ & $3.6516\times 10^{-5}(\rho=0.0000012)$ & $2.2376\times 10^{-7}$\\ 
$0.1$ & $0.02$ & $2.4573\times 10^{-6}$ & $8.6386\times 10^{-8}$ & $2.3098\times 10^{-6}(\rho=0.0000053)$ & $2.5001\times 10^{-8}$\\ 
$0.05$ & $0.01$ & $2.9692\times 10^{-7}$ & $2.6278\times 10^{-9}$ & $2.8803\times 10^{-7}(\rho=0.0000151)$ & $4.5722\times 10^{-9}$ \\ 
$0.02$ & $0.005$ & $3.6163\times 10^{-8}$ & $4.9621\times 10^{-10}$ & $3.5804\times 10^{-8}(\rho=0.0000589)$ & $6.8026\times 10^{-10}$ \\ \hline
\end{tabular}%

\section{Conclusion}

\noindent
In the study, we derive a collocation method based on exponential B-spline functions combined with strong stable Crank-Nicolson implicit method. In order to check the validity and accuracy of the method, we solve an initial boundary value problem constructed with nonlinear Klein-Gordon equation. The accuracy of the method is measured by computing the maximum error norm. Comparison with the results obtained by polynomial B-splines in \cite{boz} shows that the proposed method is capable and generates more accurate results. The relative changes of the lowest two conservation laws are also good indicators of good agreement of the numerical results with analytical ones. Different from polynomial B-splines, the change of parameter in the exponential B-splines provides recovery in the results. 

\noindent
Acknowledgements: A brief part of this study was orally presented in \textit{the International Conference on Analysis and its Applications-ICAA2016,Kırşehir, Turkey, 2016.}

\end{document}